 \documentclass[11pt,twoside]{article}
 \usepackage[dvips]{graphics}
 \usepackage{amsmath,graphicx,amsfonts,amsrefs,amssymb}
 \usepackage{setspace, enumitem, mathtools, float, yfonts}
 \usepackage{graphicx}
 \usepackage[table]{xcolor}
 \def\qed{\hfill\rule{1ex}{1ex}\\}

\begin{document}
 
 \title{Multi-scale spectral methods for bounded radially symmetric capillary surfaces}
  \author{ Jonas Haug \footnote{Department of Mathematics, Texas State University, 601 University Dr., San Marcos, TX 78666, qcd5@txstate.edu} and
 Ray Treinen \footnote{Department of Mathematics, Texas State University, 601 University Dr., San Marcos, TX 78666, rt30@txstate.edu}
 }
 \maketitle

\begin{abstract}
We consider radially symmetric capillary surfaces that are described by bounded generating curves.  We use the arc-length representation of the differential equations for these surfaces to allow for vertical points and inflection points along the generating curve.  These considerations admit capillary tubes, sessile drops, and fluids in annular tubes as well as other examples. 

We present a multi-scale pseudo-spectral method for approximating  solutions of the associated boundary value problems based on interpolation by Chebyshev polynomials.   The multi-scale approach is based on a domain decomposition with adaptive refinements within each sub-domain.
\end{abstract}


\section{Introduction}
\label{intro}

In a recent work the second author introduced an adaptive Chebyshev spectral method for  radially symmetric capillary surfaces.  These radially symmetric capillary surfaces give rise to  boundary value problems for systems of nonlinear ODEs.  Treinen \cite{Treinen2023} found approximations of the solutions to these nonlinear systems by using Newton's method to generate a sequence of linear problems  and the corresponding sequence  of solutions,  when convergent,  approach the solution of the nonlinear problem.  The solutions of the linear problems were computed with Chebyshev spectral methods.  There were checks put in place to ensure that the Newton steps were converging to a solution of the nonlinear problem within a prescribed tolerance, and if that was not achieved then the resolution of the underlying approximation scheme was adaptively increased,  repeating until convergence was achieved within the requested tolerance.
Overall the performance was fast and robust.  However, in that work there were identified some cases where the adaptive step needed a rather ungainly number of points, and the performance suffered.  There were also problems where that algorithm did not converge.   These cases are the subject of the current study.  We have implemented multi-scale approaches to these problems where the problematic regions  are isolated and treated with separate adaptive approaches.

We will be using Matlab with Chebfun \cite{Chebfun} to facilitate some aspect of our work, though we will not be using all of the automation implemented in Chebfun.  We do use Chebfun to generate our differentiation matrices, and for plotting our generating curves with barycentric interpolation in a way that puts visible points at the Chebyshev points along the curves.

We will focus on two prototype problems:
\begin{itemize}
	\item {\bf P1}
	Simply connected interfaces that are the image of a disk, and
	\item {\bf P2}
	Doubly connected interfaces that are the image of an annulus.
\end{itemize}
Within these two prototype problems we will treat subsets of problems that have been either computationally expensive to solve or where convergence altogether fails.  Before we discuss our multiscale approach to specific problems, we give some brief background on the underlying mathematical problems and summarize the algorithm from \cite{Treinen2023}.

We have released Matlab implementations of these algorithms under an open source license.  These files are available in a GitHub repository at \\
\texttt{https://github.com/raytreinen/Multiscale-Spectral-Method-Capillarity.git} \\

\section{Mathematical framework for the prototype problems}
\label{Details}

The mean curvature of a capillary surface is proportional to its height  $u$, where the height is measured over some fixed reference level.  We will restrict our attention to embedded solutions of the mean curvature equation
$$
2H = \kappa u,
$$
where $H$ is  the mean curvature of the surface, $\kappa = \rho g / \sigma > 0$ is the capillary constant with $\rho$ defined to be the difference in the densities of the fluids, $g$ the gravitational constant, and $\sigma$ the surface tension. See Finn \cite{ecs}.  This equation is sometimes known as the Young-Laplace equation, or the capillary equation.  We defer to  Finn \cite{ecs} and Treinen \cite{Treinen2023} for background,  physical interpretations of these interfaces, and further references.

A radial symmetric capillary surface can be described by a system of three nonlinear ordinary differential  equations, parametrized by the arc-length $s$:
\begin{eqnarray}
\frac{dr}{ds} &=& \cos\psi, \label{drds}\\
\frac{du}{ds} &=& \sin\psi, \label{duds}\\
\frac{d\psi}{ds} &=& \kappa u - \frac{\sin\psi}{r}, \label{dpsids}
\end{eqnarray}
where $r$ is the radius and $u$ is the height of the interface, and the inclination angle $\psi$ is measured from the corresponding generating curve which is described by $(r(s),u(s))$. 

The associated boundary values will come from the natural boundary conditions for each problem, and will correspond to prescribing the inclination angle(s) at one or more radii, and the details are described in what follows.

\begin{figure}[t]
	\centering
	\scalebox{0.35}{\includegraphics{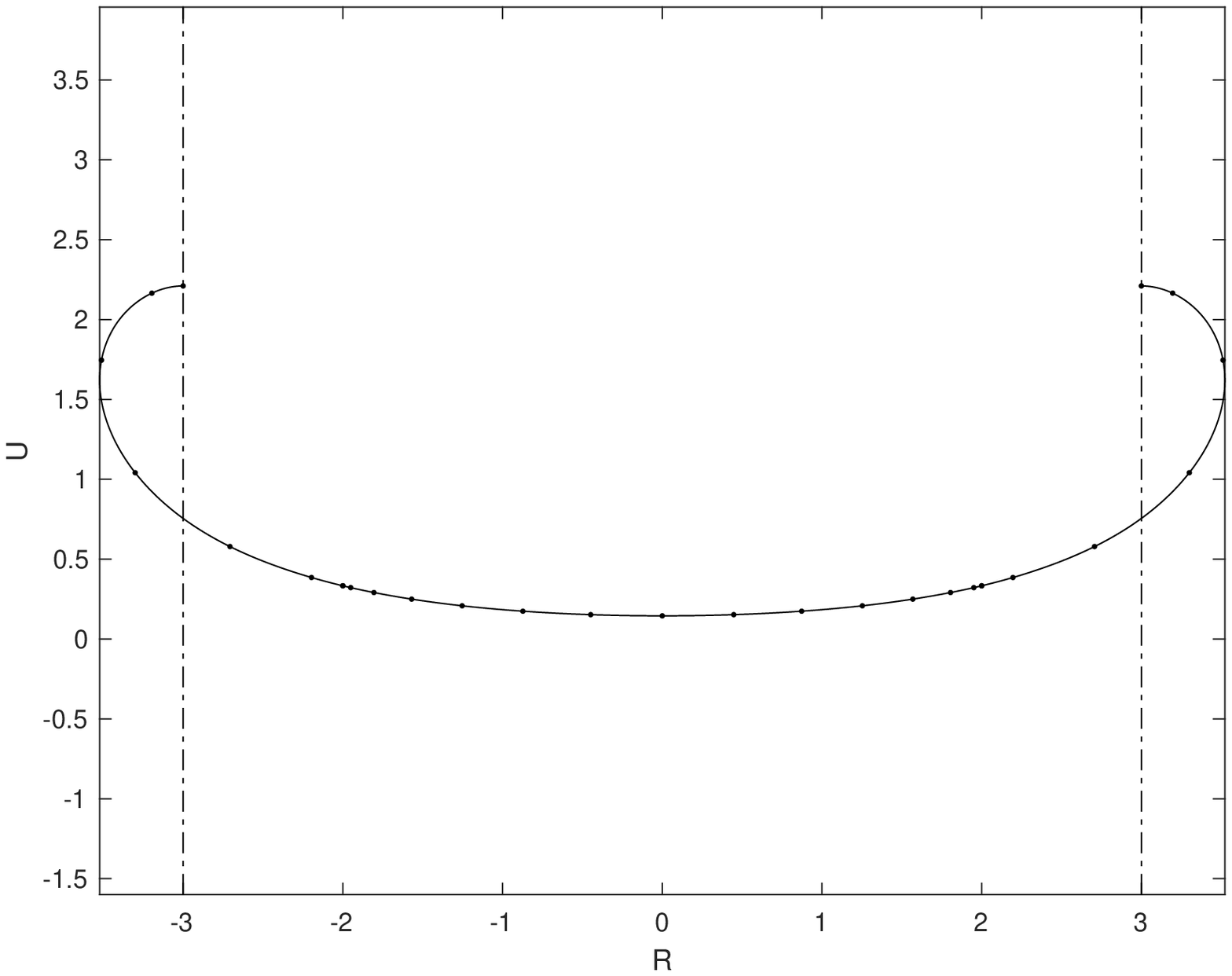}}
	\scalebox{0.35}{\includegraphics{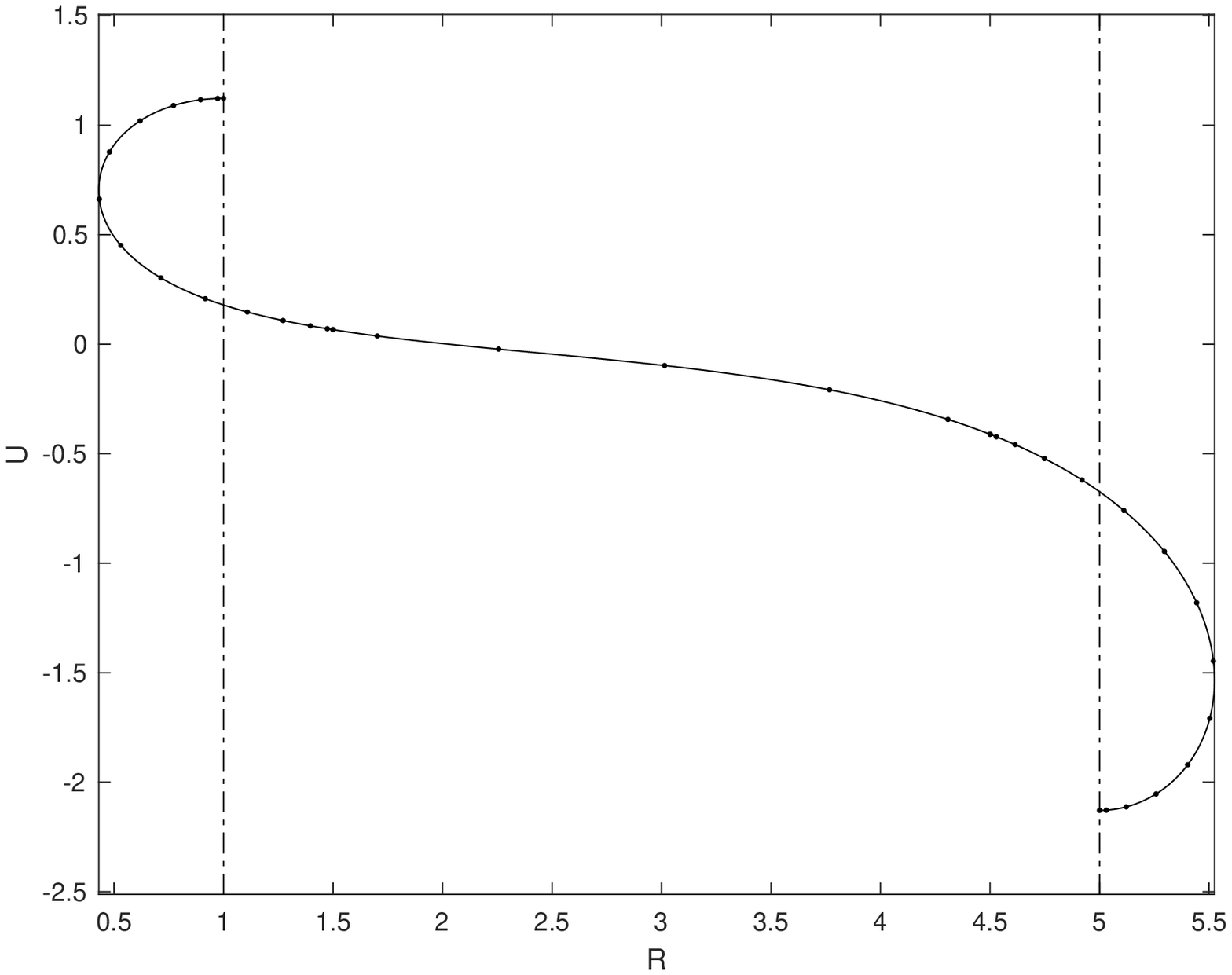}}
	\caption{A capillary surface with radius $b = 3$ and inclination angle $\psi_b = \pi$ (L) and with radii $a = 1$ and $b = 5$ with inclination angles $\psi_a, \psi_b = -\pi$ (R).  In all of our figures we include vertical lines to indicate the radius (or radii) of interest.  }
	\label{fig:intro}
\end{figure}

\subsection{ The problem {\bf P1} }
\label{subsecP1}

For radially symmetric and simply connected surfaces that are the image of a disk, we specify boundary conditions by the requirement that at some arc-length $\ell>0$ the radius $r(\ell)$ meets a prescribed value $b>0$, and the inclination angle $\psi(\ell)$ meets a prescribed value $\psi_b \in [-\pi, \pi]$.  However, the value of the arc-length $\ell$ is unknown.  We  approach this unknown parameter by rescaling the problem,  defining $\tau = s/\ell$, or $s = \ell\tau$.  Then we define
\begin{eqnarray}
R(\tau) &:=& r(\ell\tau) = r(s),  \nonumber\\
U(\tau) &:=& u(\ell\tau) = u(s), \nonumber\\
\Psi(\tau) &:=& \psi(\ell\tau) = \psi(s). \nonumber 
\end{eqnarray}
This scaling is suitable for use with Chebyshev polynomials, as the natural domains of $R,U,$ and $\Psi$ are now $[-1,1]$.  Then, using the chain rule, from   \eqref{drds}-\eqref{dpsids} we find
\begin{eqnarray}
R^\prime(\tau) - \ell\cos\Psi(\tau) &=& 0, \label{eqn:R}\\
U^\prime(\tau) - \ell\sin\Psi(\tau) &=& 0, \label{eqn:U}\\
\Psi^\prime(\tau) + \frac{\ell\sin\Psi(\tau)}{R(\tau)} - \kappa\ell U(\tau) &=& 0. \label{eqn:Psi}
\end{eqnarray}
If we define the column vector $\mathbf{v} = [R\,\, U\,\, \Psi\,  \ell]^T$, we can use \eqref{eqn:R}-\eqref{eqn:Psi} to define the nonlinear operator in the vector equation
\begin{equation}
	\tilde N(\mathbf{v}) = \mathbf{0}. \label{eqn:tN=0}
\end{equation}
We then use the boundary conditions
\begin{eqnarray}
R(1) - b &=& 0, \label{eqn:bcRb} \\
R(-1) + b &=& 0, \label{eqn:bcRmb} \\
\Psi(1) - \psi_b &=& 0, \label{eqn:bcPsib} \\
\Psi(-1) + \psi_b &=& 0, \label{eqn:bcmPsib}
\end{eqnarray}
so that we have some form of a two-point boundary value problem.  Here we are prescribing the natural boundary conditions on a cylinder of radius $b$.  See Figure~\ref{fig:intro} (L).  We also allow for the strictly parametric solution where $\pi/2 < \psi_b$, and this then provides a model for the sessile drop wetting a disc of radius $b$ under appropriate rigid transformations.  We append \eqref{eqn:tN=0} with these  boundary conditions to form the system 
\begin{equation}
{N}(\mathbf{v}) = \mathbf{0}.
\end{equation}
The resulting solution $(R,U) = (r,u)$ is not a generating curve of the surface, but it is a section.  The generating curve is given by $(R(\tau),U(\tau))$ for $0\leq\tau\leq 1$, and the generating curve has total arc-length $\ell$.  One final remark on the general nature of this problem is that \eqref{dpsids} and \eqref{eqn:Psi} have singularities at $r = R = 0$.  It is known that this singularity is removable, and the solution is analytic \cite{ecs}.  We will discuss our numerical approach to this singularity when we return to the results for this problem.  For the time being, we present an alternative form of \eqref{eqn:Psi} that is suitable for reducing the numerical error from rounding, as exaggerated by dividing by a number close to zero.  We multiply by $R$ to get
\begin{equation}
R(\tau)\Psi^\prime(\tau) + \ell\sin\Psi(\tau) - \kappa\ell R(\tau) U(\tau) = 0, \label{eqn:Psialt}
\end{equation}
which we will use when $b$ is relatively small.  When we use \eqref{eqn:Psialt} in place of \eqref{eqn:Psi} we denote the changed $N$ by $N_1$ (and $F$ and $L$ similarly below) when specifying the difference is important, and we will generically refer to those objects without subscripts when no confusion is expected.

We will later approach this nonlinear problem with a Newton method, and we will need to use the Fr\'{e}chet derivative 
$$
F(\mathbf{v}) = \frac{d N}{d\mathbf{v}}(\mathbf{v}).
$$
Given that $\mathbf{v}$ has several components, and since some of the computations in $F(\mathbf{v})$ involve derivatives with respect to $\tau$, we introduce the differential operator 
$$
D = \frac{d\, }{d\tau}, 
$$
which is applied in a block fashion to $\mathbf{v}$ so that $R^\prime(\tau) = [D\,\, 0\,\, 0\,\, 0]\mathbf{v}$, for example.  We will also have need to use an operator version of function evaluation.  We denote $D^0_\tau$ to be this operator, so that $D^0_1R = R(1)$.  With this in hand, we compute
\begin{equation} \label{eqn:N=0}
F(\mathbf{v} ) = 
\begin{bmatrix}
	D & 0 & \ell\sin\Psi & -\cos\Psi \\
	0 & D & -\ell\cos\Psi & -\sin\Psi \\
	\frac{-\ell\sin\Psi }{R^2} & -\kappa\ell & D + \frac{\ell\cos\Psi }{R} & \frac{\sin\Psi}{R} - \kappa U \\
	D^0_{-1} & 0 & 0 & 0 \\
	D^0_{1} & 0 & 0 & 0 \\
	0 & 0 & D^0_{-1} &  0 \\
	0 & 0 & D^0_{1} &  0 
\end{bmatrix}
\mathbf{v}.
\end{equation}
We will have need to solve linear systems based on the definition $F(\mathbf{v}) := L\mathbf{v}$.  For $F_1$ and $L_1$, we merely change the third row to 
\begin{equation}
	\label{multR}
\begin{bmatrix}
D\Psi - \kappa\ell U & -\kappa\ell R & RD + \ell\cos\Psi & \sin\Psi - \kappa UR
\end{bmatrix}.
\end{equation}

\subsection{ The problem {\bf P2} }
\label{subsecP2}

If we replace the symmetric and simply connected surfaces that are the image of the disk considered in Subsection~\ref{subsecP1} with doubly connected interfaces that are the image of an annulus, the scaling argument is preserved.  There are differences though, as the arc-length $s=0$ does not correspond to the radius being zero.  We reuse \eqref{eqn:R}-\eqref{eqn:Psi} and \eqref{eqn:N=0}, but we have boundary conditions
\begin{eqnarray}
	R(1) - b &=& 0, \label{eqn:bcRb2} \\
	R(-1) - a &=& 0, \label{eqn:bcRa} \\
	\Psi(1) - \psi_b &=& 0, \label{eqn:bcPsib2} \\
	\Psi(-1) - \psi_a &=& 0, \label{eqn:bcPsia}
\end{eqnarray}
where we introduce $a\in(0,b)$ and $\psi_a\in[-\pi,\pi]$.  See Figure~\ref{fig:intro} (R).  This changes the definition of $N$, but not that of $F$.  Notice that in this case $R \ne 0$, so there is no singularity.  If $|\psi_a|, |\psi_b| \leq \pi/2$, then the resulting solution curves are graphs and we have a so-called annular capillary surface in an annular tube, with the inner wall at the radius $a$ and the outer wall at the radius $b$ \cite{EKT2004}.  If this restriction of the inclination angles is removed, then the solution curve can pass through one or both of the ``walls''  \cite{Treinen2012}.   In either case, these solution curves form generating curves with total arc-length $2\ell$.

\section{The spectral method from \cite{Treinen2023}}
\label{Cheby}

We summarize the spectral method developed by Treinen \cite{Treinen2023}, as it is the fundamental building block of our approach.
We also point out that this approach has its roots in Trefethen's monograph \cite{Trefethen2000}, Driscoll and Hale \cite{DriscollHale2016}, and  Aurentz and Trefethen \cite{AurentzTrefethen2017}, and if the reader is not yet familiar with these works, he or she is strongly encouraged to read these first, perhaps starting with \cite{AurentzTrefethen2017}.  

\begin{figure}[!h]
	\centering
	\scalebox{0.37}{\includegraphics{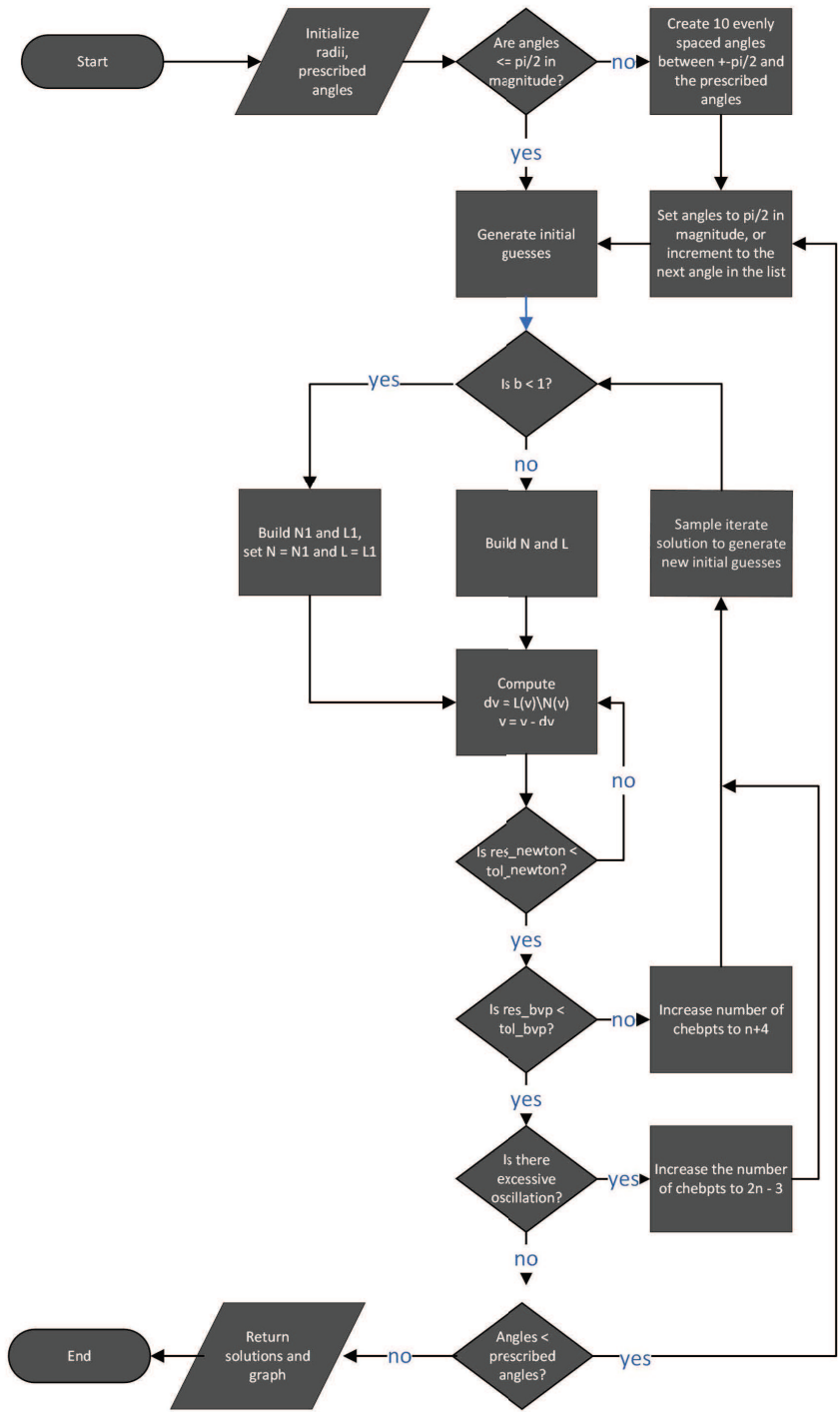}}
	\caption{The flowchart describing the general algorithm.}
	\label{flowchart}
\end{figure}

Here we present the basic ideas with the backdrop of $\mathbf{P1}$.   The flowchart in Figure~\ref{flowchart} gives the framework of the algorithm, and applies to  both $\mathbf{P1}$ and $\mathbf{P2}$.  Here we will merely discuss the key points of this algorithm.

We will treat the nonlinearity of $N$ with a Newton method.   The core of the algorithm is to construct $N$ and $L$, and then we solve the linear equation $L(\mathbf{v}) \, d\mathbf{v} = - N(\mathbf{v})$ for  $d\mathbf{v}$ to build the update $\mathbf{v}_{\mathtt{new}} = \mathbf{v} + d\mathbf{v}$.  This iterative process continues until some convergence criteria is met.  This process needs an initial guess to begin, and if that initial guess is sufficiently close to the solution, then the convergence is known to be at least quadratic.

We will  discuss the initial guesses for Newton's method when we consider specific problems below.   However, we will  not discuss the barriers constructed  in \cite{Treinen2023} which were used to enforce convergence to physically correct configurations.  These barriers are still used here.

We will need Chebyshev differentiation matrices in what follows.  These matrices can be realized as representing the linear transformation between two vectors of data corresponding to particular grid points points, say $\mathbf{f}$ being mapped to $\mathbf{f^\prime}$.  The data at the grid points correspond to the interpolating polynomials for a function $f$ and its derivative $f^\prime$, where the data is sampled at Chebyshev grid points $x_j = \cos(\theta_j)\in[-1,1]$ with angles $\theta_j$ equally spaced angles over $[0,\pi]$.  Of course, these grid points do not need to be fixed, as multiple samplings can be used.  

The basic building blocks of  both the nonlinear equation $N(\mathbf{v}) = \mathbf{0}$ and the linearized equation $L d\mathbf{v} = -N$  are based on $D^0$ and $D$, which we implement using the Chebfun commands
\begin{verbatim}
D0 = diffmat([n-1 n],0,X);
D1 = diffmat([n-1 n],1,X);
\end{verbatim}
where $X = [-1,1]$,  $n$ is the number of Chebyshev points we are using, and the input 0 or 1 indicates the number of derivatives.  Since $\texttt{D0}$  is rectangular, it becomes a $(n-1)\times n$ identity matrix interpreted as a dense ``spectral down-sampling'' matrix implemented as interpolating on an $n$-point grid followed by sampling on an $(n-1)$-point grid.    We  sparsely construct our operators $N$ and $L$ using these components in a block fashion.
These building blocks are found by multiplying a linear operator times the vector $\mathbf{v}$ that contains the solution components.  Then we use these building blocks to put together $N$ and $L$.  If $b$ is small, then we construct $N_1$ and $L_1$ where we have no excessive numerical error due to rounding near the singularity.

We always take care that the number of Chebyshev points are chosen so that there is no evaluation of  $r = 0$.  This removable singularity is avoided in exactly the same way Trefethen described in \cite{Trefethen2000} when discussing the radial form of Laplace's equation.

The basic loop is
\begin{verbatim}
	while res_newton > tol_newton
    dv = L(v)\N(v);
	    v = v - dv;
	    res_newton = norm(dv,'fro')/norm(v,'fro');
	end
\end{verbatim}
The tolerance $\texttt{tol\_newton}$ must be met, and we use the relative error measured by the Frobeneous norm.  For all of the examples in this paper, we used $\texttt{tol\_newton = 1e-13}$.

If the Newton's method fails to converge within a specified maximum number of iterations, we then increase the number of Chebyshev points and we sample the current state of the iteration onto the new Chebyshev points, reinitialize the operators $N$ and $L$ as above, and we enter another loop for Newton's method.  Further, we include this process in an outer  loop that also tests the relative error of the iterates even if the Newton's method converges.  We use the Frobenious norm to compute $\texttt{res\_bvp}$ as $||N(v)||/||v||$.  If this residual is greater than the prescribed tolerance, then we also increase the number of Chebyshev points.  We used a prescribed tolerance of $\texttt{tol\_bvp = 1e-12}$ in all of the examples in this paper.  It may also be that there is excessive polynomial oscillation across the Chebyshev points, and there are not enough of these points to accurately resolve the solution,  then we again increase the number of Chebyshev points, however this increase is more aggressive.

If the boundary conditions are specified so that the solution curve will not be a graph over a base-domain, such as when $\psi_b > \pi/2$, then we use a method of continuation.  This applies to most of the problems considered in this paper.  We harvest the sign of $\psi_b$ and in place of  that boundary condition, we solve ten problems with linearly spaced boundary angles between $\pi/2$ and $|\psi_b|$ using the  converged solution at each step as an initial guess for the next step.  We preserve the number of needed Chebyshev points from the adaptive algorithm while moving from one step to the next.  
In $\mathbf{P2}$, if needed, we do this for the inclination angles at $a$ and $b$ simultaneously.

\section{A multi-scale approach via domain decomposition}
\label{multiscale}

We  partition the arc-length domains into two or three sub-domains, depending on the problem.  The general approach is to identify a parameter that marks the boundary between the sub-domains and then conditions are put in place so that the solutions on the sub-domains extend across this artificial boundary as a solution for the global problem.

For many of the problems  we choose the parameter to be a radius.  We use the setting of {\bf P1} to describe the details of the process for that choice of parameter.  For some problems we find the parameter of the inclination angle to be a better choice for marking the boundary of the sub-domains, and we use the setting of {\bf P2} to describe the process for that choice of parameter.

In either case, we build a vector $\mathbf{v}$, and construct appropriate operators $N(\mathbf{v})$ and $F(\mathbf{v}) = L(\mathbf{v})\mathbf{v}$ for the partitioned domains with the matching conditions included.  Then we follow the basic algorithm as described in the flowchart in Figure~\ref{flowchart} where the relative errors and adaptive steps are taken corresponding to each sub-domain.  

It should be noted that while this is simple enough to describe, the computer code needed to implement these concepts becomes somewhat lengthy and care must be taken in the adaptive procedures for generating $\mathbf{v}$, $N$, and $L$.  We will provide more detail in what follows.

\subsection{{\bf P1} Study}
\label{Study1}

For this case we require that $\psi = \psi_b$ at $r = b$, and, by symmetry, $\psi = -\psi_b$ at $r = -b$.  The base code \cite{Treinen2023} works well in general, however we are able to improve the performance if $\psi_b> \pi/2$ and $b \gg 1$.  We then choose a $\delta > 0$ and pick the parameters to mark the sub-domain boundaries at $r = \delta -b$ and $r = b-\delta$.  This partitions the arc-length domain into three sub-domains.  Denote those domains by $\Omega_1$, $\Omega_2$, and $\Omega_3$ respectively as one moves  from left to right along the solution curve.   We will also use these subscripts to denote the restrictions of solutions $(r,u,\psi)$ on those sub-domains.  The arc-lengths at the boundary points are denoted by $\ell_1$ and $\ell_2$.   

Physically the height $u$ and inclination angle $\psi$ on $\Omega_2$ will both be approximately zero, and the radius will be approximately measured by a translation of the arc-length.  Thus relatively very few Chebyshev points are needed to resolve the solution to within the tolerance requested.  On both $\Omega_1$ and $\Omega_3$ this no longer is the case, and the adaptive method of continuity described in Section~\ref{Cheby} is used to resolve the solutions there.  The conditions over the boundaries are
\begin{eqnarray}
	u_1(\ell_1) = u_2(\ell_1) , \\
	\psi_1(\ell_1) = \psi_2(\ell_1),
\end{eqnarray}
and
\begin{eqnarray}
	u_2(\ell_2) = u_3(\ell_2),  \\
	\psi_2(\ell_2) = \psi_3(\ell_2).
\end{eqnarray}

Some explanation of initial guesses for the underlying Newton's method is in order, though we will be brief as \cite{Treinen2023} contains the formulas for the functions we need.  On $\Omega_1$ we will use a circular arc that meets the boundary conditions at $r = -b$ and is horizontal at $r = -b + \delta$.  On $\Omega_2$ we use a line with zero height.  On $\Omega_3$ we again use a circular arc that meets the boundary conditions at $r = b$ and is horizontal at $r = b - \delta$.  We also prescribe the heights of these circular arcs so that the heights match at the boundaries of the sub-domains.  Since we are primarily interested in problems where $|\psi_b| > \pi/2$, we use these guesses to produce a converged solution for $|\psi_b| = \pi/2$ and then we use the method of continuation to feed that converged solution to problems with a sequence of increasing boundary angles up to $|\psi_b|$.

The computational vector is 
$\mathbf{v} = [R_1\, R_2 \, R_3\, U_1 \, U_2 \, U_3\, \Psi_1\,  \Psi_2\, \Psi_3\, \ell_1\, \ell_2\, \ell_3]^T$
and the mathematical description of $N$ is given by
\begin{equation}
	\label{bigN}
	N(\mathbf{v}) = 
\begin{bmatrix}
	R_1^\prime(\tau_1) - \ell_1\cos\Psi_1(\tau_1)\\
		R_2^\prime(\tau_2) - \ell_2\cos\Psi_2(\tau_2)\\
			R_3^\prime(\tau_3) - \ell_3\cos\Psi_3(\tau_3)\\
	U_1^\prime(\tau_1) - \ell_1\sin\Psi_1(\tau_1)\\
		U_2^\prime(\tau_2) - \ell_2\sin\Psi_2(\tau_2)\\
			U_3^\prime(\tau_3) - \ell_3\sin\Psi_3(\tau_3)\\
	\Psi_1^\prime(\tau_1) + \frac{\ell_1\sin\Psi_1(\tau_1)}{R_1(\tau_1)} - \kappa\ell_1 U_1(\tau_1) \\
		\Psi_2^\prime(\tau_2) + \frac{\ell_2\sin\Psi_2(\tau_2)}{R_2(\tau_2)} - \kappa\ell_2 U_2(\tau_2) \\
			\Psi_3^\prime(\tau_3) + \frac{\ell_3\sin\Psi_3(\tau_3)}{R_3(\tau_3)} - \kappa\ell_3 U_3(\tau_3) \\
	R_1(-1) + b \\
		R_1(1) + b - \delta\\
	\Psi_1(-1) + \psi_b  \\
	\Psi_1(1) - 	\Psi_2(-1) \\
			R_2(-1) + b - \delta\\
				U_1(1) - U_2(-1) \\
			R_2(1) - b + \delta\\
				\Psi_2(1) - 	\Psi_3(-1) \\
			R_3(-1) - b + \delta\\
			U_2(1) - U_3(-1) \\
R_3(1) - b \\
\Psi_3(1) - \psi_b 
\end{bmatrix}.
\end{equation}
For $i=1,2,3$ we allocate in each of these sub-domains an independent variable $\tau_i$,   and at the appropriate step of the algorithm the number of gridpoints in each $\tau_i$ will be adaptively increased if the tolerances for the residual of $N$ restricted to $\Omega_i$ is not met, as indicated in Figure~\ref{flowchart}.  We scale the problem as in Section~\ref{Cheby} to determine these arc-lengths $\ell_i$, for $i = 1,2,3$ in the process of solving the problem. 
Then $L$ is given by the block matrix
\begin{equation}
\label{kronL}
L = 
\begin{bmatrix}
M_{11} & M_{12}  \\
M_{21} & M_{22}  \\
M_{31} & M_{32}  \\
\end{bmatrix},
\end{equation}
where $M_{11} = I\otimes D$ is given using the Kronecker product where $I$ is the identity  on $\mathbb{R}^{6\times 6}$, 
$$
M_{12}  = 
\begin{bmatrix}
\ell_1\sin\Psi_1 & 0 & 0&-\cos\Psi_1 & 0 & 0  \\
0 & \ell_2\sin\Psi_2  & 0& 0& -\cos\Psi_2  & 0  \\		
0 & 0 & \ell_3\sin\Psi_3 & 0 & 0&-\cos\Psi_3  \\
-\ell_1\cos\Psi_1 & 0& 0& -\sin\Psi_1& 0& 0 \\
0 & -\ell_2\cos\Psi_2 & 0& 0& -\sin\Psi_2& 0 \\		
0 & 0 &  -\ell_3\cos\Psi_3 & 0& 0& -\sin\Psi_3 
\end{bmatrix},
$$

$$
M_{21}  = 
\begin{bmatrix}
\frac{-\ell_1\sin\Psi_1 }{R_1^2}& 0& 0 & -\kappa\ell_1 & 0& 0 \\
0 & \frac{-\ell_2\sin\Psi_2 }{R_2^2}& 0& 0 & -\kappa\ell_2 & 0 \\
0 & 0 & \frac{-\ell_3\sin\Psi_3 }{R_3^2}& 0& 0 & -\kappa\ell_3
\end{bmatrix},
$$
$$
M_{22}  = 
\begin{bmatrix}
 D  + \frac{\ell_1\cos\Psi_1 }{R_1}& 0& 0&  - \kappa U_1 & 0& 0\\
 0& D  + \frac{\ell_2\cos\Psi_2 }{R_2}& 0& 0&  - \kappa U_2 &  0\\
 0& 0& D  + \frac{\ell_3\cos\Psi_3 }{R_3}& 0& 0&  - \kappa U_3 
\end{bmatrix},
$$
and the blocks corresponding to the boundary conditions and matching conditions are given by
$$
M_{31}  = 
\begin{bmatrix}
D^0_{-1} & 0 & 0 & 0 & 0& 0 \\
D^0_{1} & 0 & 0 & 0 & 0& 0 \\
0 & 0 & 0 & 0 & 0& 0 \\
0 & 0 & 0 & 0 & 0& 0 \\
0 & D^0_{-1} & 0 & 0 & 0 & 0 \\
0 & 0 & 0 & D^0_1 & -D^0_{-1} & 0 \\
0 & D^0_{1} & 0 & 0 & 0 & 0 \\
0 & 0 & 0 & 0 & 0 & 0 \\
0 & 0 & D^0_{-1} & 0 & 0 & 0  \\
0 & 0 & 0 & 0 & D^0_1 & -D^0_{-1} \\
0 & 0 & D^0_{1} & 0 & 0 & 0 \\
0 & 0 & 0 & 0 & 0 & 0 
\end{bmatrix},
$$
and
$$
M_{32}  = 
\begin{bmatrix}
0 & 0 & 0 & 0 & 0& 0 \\
0 & 0 & 0 & 0 & 0& 0 \\
D^0_{-1} & 0 & 0 & 0 & 0& 0 \\
D^0_{1} & -D^0_{-1} & 0 & 0 & 0& 0 \\
0 & 0 & 0 & 0 & 0& 0 \\
0 & 0 & 0 & 0 & 0& 0 \\
0 & 0 & 0 & 0 & 0& 0 \\
0 & D^0_{1} & -D^0_{-1} & 0 & 0 & 0 \\
0 & 0 & 0 & 0 & 0& 0 \\
0 & 0 & 0 & 0 & 0& 0 \\
0 & 0 & 0 & 0 & 0& 0 \\
0 & 0 & D^0_1 & 0 & 0 & 0
\end{bmatrix}.
$$
With these objects in hand we  use the general algorithm from the flowchart in Figure~\ref{flowchart}.

\begin{figure}[t]
	\centering
	\scalebox{0.35}{\includegraphics{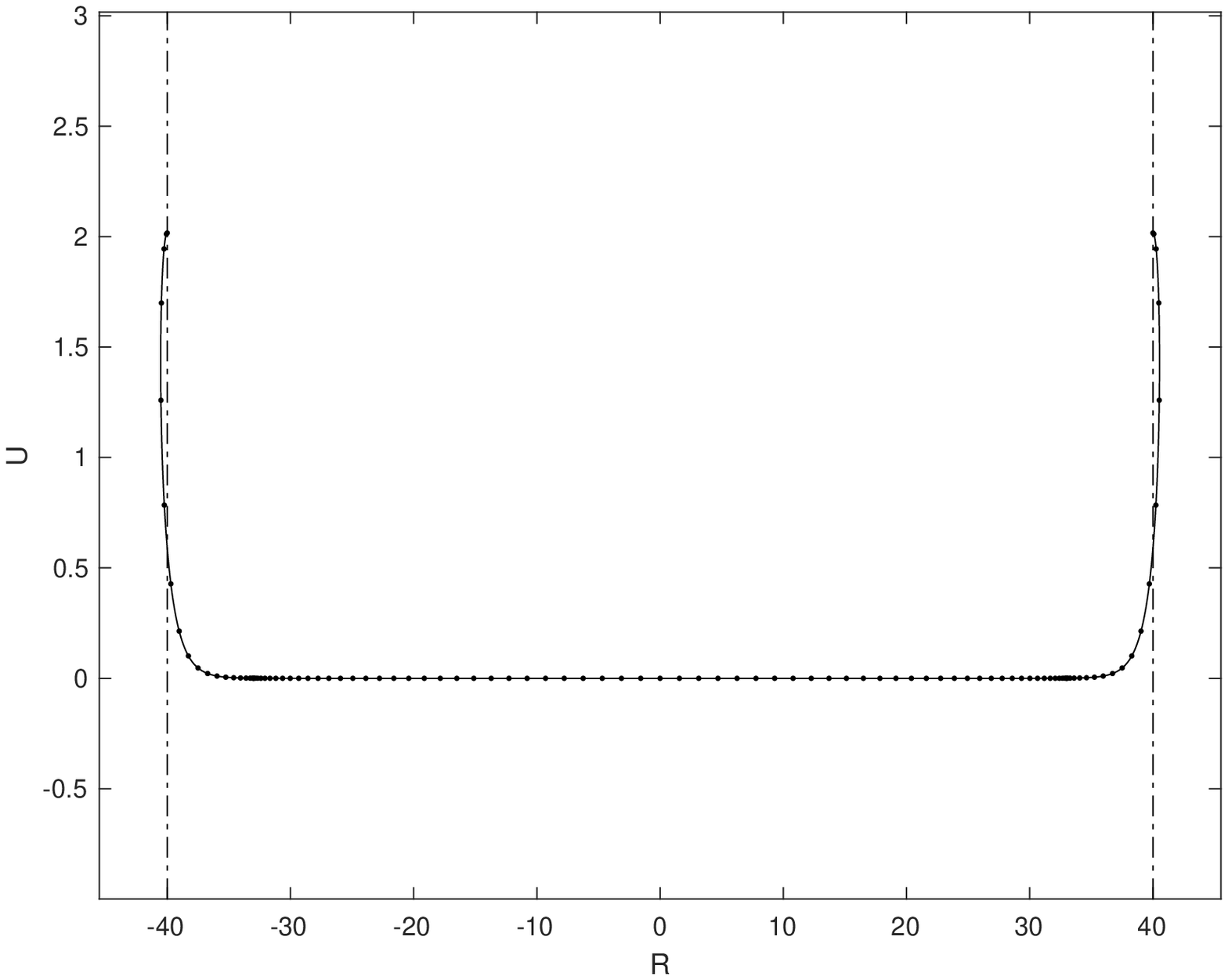}}
	\scalebox{0.35}{\includegraphics{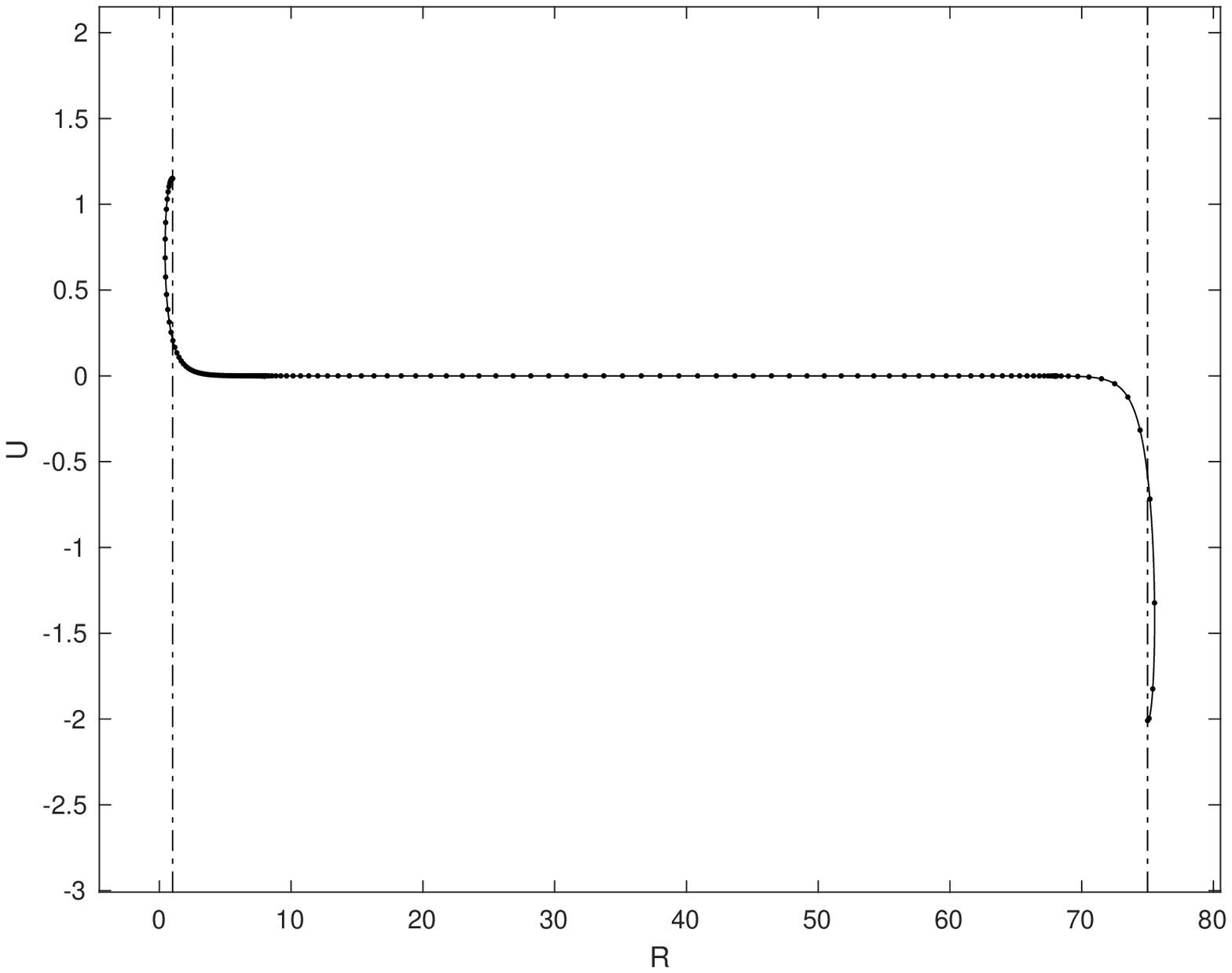}}
	\caption{A capillary surface with radius $b = 40$ and inclination angle $\psi_b = \pi$ (L) and with radii $a = 1$ and $b = 75$ with inclination angles $\psi_a, \psi_b = -\pi$ (R).  The vertical axes have been scaled to highlight these examples.}
	\label{fig:long}
\end{figure}

For measures of performance, we are primarily interested in the size of the linear systems used and the number of Newton iterations needed for convergence.  To give an approximation of the sizes of these problems, denote $n_v$ to be the length of the computational vector $\mathbf{v}$, giving the size of the linear system.  We denote the number of Newton iterations by $n_N$, and this roughly indicates how many times the algorithm solves a system of size $n_v\times n_v$.  Note that since the size of $\mathbf{v}$ is potentially growing in the adaptive phases, these numbers only give an upper bound on the complexity of the problem.  While it would be possible to compute the exact numbers, we find that this approach illuminates the process well enough without the extra detail.   Of course we are also interested in  the configurations where the multi-scale algorithm  converges when the program from \cite{Treinen2023} fails to meet the same tolerances.  Some examples of these experiments are collected in Table~\ref{table1}.  Figure~\ref{fig:long} (L) shows the last example from the table.  We focused on problems with large radii $b$, so the approach indicated by \eqref{multR} was not used.

\begin{table}[t]
	\centering
	\begin{tabular}{||l|l|l|l|l|l|l||}
		\cline{1-7}
		$b$ & $\psi_b$ & base code $n_v$ & base code $n_N$ & MS code $n_v$ & MS code $n_N$& $\delta$ \\
		\hline \hline
		11 & $\pi$ & 154 & 53 & 90 & 130 & 2 \\
		29 & $3\pi/8$ & 298 & 16 & 246 & 132 & 3 \\
		29 & $\pi$ & 586 & 58 & 246 & 132 & 3 \\
		40 & $\pi$ & dnc & & 270 & 139  & 5 \\
		\hline
	\end{tabular}
	\caption{Solutions that are the image of a disk with radius $b$ are collected for a selection of inclination angles. The size of the corresponding linear system is given by  $n_v$, for the length of the computational vector $\mathbf{v}$, while the number of Newton iterations is given by $n_N$.  We label the columns corresponding to the multi-scale code by MS.  The value of $\delta$ used for each experiment is included.  If the algorithm did not converge it is labeled as dnc.}
	\label{table1}
\end{table}

\subsection{{\bf P2} Study}
\label{Study2}

For the class of problems {\bf P2} that are the image of an annulus, given $0<a<b<\infty$, we  require that $\psi = \psi_a$ at $r = a$ and $\psi = \psi_b$ at $r = b$.   We will use our approach on two types of configurations.  First where $a \ll \infty$ and $b$ is somewhat large, and second where $0<a\ll 1$ and $b$ is moderately sized.  To accomplish this we will use three different domain decompositions.

We consider first a case with three radial zones or sub-domains, and refer to it by 3RZ below.  We choose a $\delta > 0$ and pick our parameters to mark the boundaries at $r = a + \delta$ and $r = b-\delta$ to partition the arc-length domain into three sub-domains just as we did in Section~\ref{Study1}.  
The primary difference in this case is that the last 12 rows of $N(\mathbf{v})$ corresponding to the boundary conditions and matching conditions are replaced with 
$$
\begin{bmatrix}
R_1(-1) -a \\
R_1(1)  - a + \delta\\
\Psi_1(-1) - \psi_a  \\
\Psi_1(1) - 	\Psi_2(-1) \\
R_2(-1) -a + \delta\\
U_1(1) - U_2(-1) \\
R_2(1) - b + \delta\\
\Psi_2(1) - 	\Psi_3(-1) \\
R_3(-1) - b + \delta\\
U_2(1) - U_3(-1) \\
R_3(1) - b \\
\Psi_3(1) - \psi_b 
\end{bmatrix}.
$$
and $L$ is unchanged.  If $a<1$ then the row replacement of \eqref{multR} is used in $\Omega_1$, and we will illustrate those details in the following cases.   Table~\ref{table2} and Table~\ref{table2pi} include some experiments that highlight the performance here.  Table~\ref{table2} has cases with $\psi_b<0$ and Table~\ref{table2pi} has cases with $\psi_b > 0$.  Figure~\ref{fig:long} (R) shows the last example from Table~\ref{table2}.

\begin{table}[t]
	\centering
	\begin{tabular}{||l|l|l|l|l|l|l|l|l|l|l||}
		\cline{1-11}
		$a$ & $b$ & $\psi_a$ & $\psi_b$ & base  $n_v$ & base  $n_N$ & 3RZ  $n_v$& 3RZ $n_N$& 2RZ $n_v$& 2RZ $n_N$ &$\delta$  \\
		\hline \hline
		1 & 5   & $-\pi/2$         & $-\pi/2$         & 46      & 8   & 114 & 69   & 80    & 6   & 0.4 \\		

		1 & 5   & $-15\pi/16$  & $-15\pi/16$ & 46      & 58   & 114 & 175   & 80    & 117   & 0.4 \\		
		1 & 5   & $-\pi$         & $-\pi$         & 82      & 63   & 114 & 189   & 80    & 138   & 0.4 \\		

		1 & 44 & $-\pi/2$        & $-\pi/2$         & 586    & 24 & 198 & 67   & 248 & 21 & 3 \\		
		1 & 44 & $-15\pi/16$ & $-15\pi/16$ & 1162 & 71 & 198 & 158 & 476 & 134 & 3 \\		
		1 & 44 & $-\pi$             & $-\pi$             & 1174 & 73 & 198 & 171 & 464& 152 & 3 \\		
		1 & 47 & $-\pi/2$        & $-\pi/2$          & dnc    &       & 198 & 69   & dnc &           & 5 \\		
		1 & 47 & $-\pi$            & $-\pi$              & dnc     &       & 246& 183 &  dnc &          & 5 \\		
		1 & 75 & $-\pi$            & $-\pi$              & dnc     &       & 462& 186 &  dnc &          & 7 \\		
		\hline
\end{tabular}
\caption{Solutions that are the image of an annulus with radii $a$ and $b$, a selection of inclination angles focused on $\psi_b = -\pi$, and the total number of points as before.}
\label{table2}
\end{table}

\begin{table}[t]
	\centering
	\begin{tabular}{||l|l|l|l|l|l|l|l|l|l|l||}
		\cline{1-11}
		$a$ & $b$ & $\psi_a$ & $\psi_b$ & base  $n_v$ & base  $n_N$ & 3RZ  $n_v$& 3RZ $n_N$& 2RZ $n_v$& 2RZ $n_N$ &$\delta$  \\
\hline\hline
		1	& 5	& - $\pi$	& $\pi$ & 58* & 56* & 114 & 151 & 80	& 135 & 0.4\\
						
		1	& 18& - $\pi$	& $\pi$ & 298     & 65 & 126 & 157 & 152	& 137 & 2\\
								
		1	& 44& - $\pi$	& $\pi$ & 1174 & 808  & 210 & 174 & dnc	&           & 3 \\	
		1	& 47 & - $\pi$	& $\pi$ & dnc    &          & 318 & 179 & dnc	&          & 3\\
		1	& 75& - $\pi$	& $\pi$ & dnc    &          & 510 & 177 & dnc	&         & 5\\
		\hline
	\end{tabular}
	\caption{Solutions that are the image of an annulus with radii $a$ and $b$, a selection of inclination angles focused on $\psi_b = \pi$, and the total number of points as before.  For the base code example marked with *, we increased the initial number of Chebyshev points $n$ from 15 to 19 to eliminate excessive oscillation. }
	\label{table2pi}
\end{table}

For next case under consideration, we split the domain into two regions as described in Section~\ref{multiscale} where we again use the radius as the parameter that marks the boundary between the two sub-domains.  We use $\Omega_1$ for the arc-lengths on the left of $r = a + \delta$, and $\Omega_2$ for the arc-lengths to the right of $a + \delta$.  We will refer to this case by 2RZ for two radial zones.  Here we are primarily interested in small values of $a$ and modest values of $b$, as we can see from the performance indicated in Table~\ref{table2}.     We have 
$\mathbf{v} = [R_1\, R_2 \, U_1 \, U_2 \, \Psi_1\,  \Psi_2\,  \ell_1\, \ell_2]^T$ and 
\begin{equation}
	\label{bigN}
	N(\mathbf{v}) = 
	\begin{bmatrix}
		R_1^\prime(\tau_1) - \ell_1\cos\Psi_1(\tau_1)\\
		R_2^\prime(\tau_2) - \ell_2\cos\Psi_2(\tau_2)\\
		U_1^\prime(\tau_1) - \ell_1\sin\Psi_1(\tau_1)\\
		U_2^\prime(\tau_2) - \ell_2\sin\Psi_2(\tau_2)\\
		R_1(\tau_1)\Psi_1^\prime(\tau_1) + \ell_1\sin\Psi_1(\tau_1) - \kappa\ell_1 R_1(\tau_1) U_1(\tau_1) \\
		\Psi_2^\prime(\tau_2) + \frac{\ell_2\sin\Psi_2(\tau_2)}{R_2(\tau_2)} - \kappa\ell_2 U_2(\tau_2) \\
		R_1(-1) -a \\
		R_1(1)  - a - \delta\\
		\Psi_1(-1) + \psi_b  \\
		\Psi_1(1) - 	\Psi_2(-1) \\
		R_2(1) - b \\
		\Psi_2(1) - \psi_b 
	\end{bmatrix}.
\end{equation}
Then we compute the Fr\'{e}chet derivative $F(\mathbf{v}) = L\mathbf{v}$ to find the following expression for $L$:
$$
\begin{bmatrix}
	D & 0 & 0 & 0 & \ell_1\sin\Psi_1 & 0 & -\cos\Psi_1 & 0 \\
	0 &	D & 0 & 0 & 0 & \ell_2\sin\Psi_2 & 0 & -\cos\Psi_2  \\
	0 & 0 &	D  & 0 & -\ell_1\cos\Psi_1 & 0 & -\sin\Psi_1 & 0 \\
	0 & 0 &	0  & D & 0 & -\ell_2\cos\Psi_2 & 0 & -\sin\Psi_2 \\
	\Psi_1^\prime - \kappa\ell_1 U_1& 0 & -\kappa\ell_1 R_1 & 0 & D + \ell_1\cos\Psi_1 & 0 & \sin\Psi_1 - \kappa R_1U_1 & 0 \\
	0 & \frac{-\ell_2\sin\Psi_2 }{R_2^2}& 0 & -\kappa\ell_1 & 0&  D + \frac{\ell_2\cos\Psi_2 }{R_2} & 0 & \kappa U_2 \\
	D^0_{-1} & 0 & 0 & 0 & 0 & 0 & 0 & 0 \\
	D^0_{1} & 0 & 0 & 0 & 0 & 0 & 0 & 0 \\
    0 & 0 & 0 & 0 & D^0_{-1} & 0 & 0 &  0 \\
    0 & 0 & 0 & 0 & D^0_{1} & -D^0_{-1} & 0 &  0 \\
    0  & D^0_{1} & 0 & 0 & 0 & 0 & 0 & 0 \\
        0 & 0 & 0 & 0 & D^0_{1} & 0 & 0 &  0 
\end{bmatrix}.
$$

The initial guesses for the underlying Newton's method are somewhat different in this case and the next.  For the problems considered up to now we used knowledge of the expected geometry of the solution to guide the initial guesses.  In the problems with small $a$ and modest $b$ where $|\psi_a|,|\psi_b| > \pi/2$ we do not have particularly helpful estimates for the solution.  So we use the base code from \cite{Treinen2023} to generate a solution where the angles are $sgn(\psi_a)\pi/2$ at $r=a$ and $sgn(\psi_b)\pi/2$ at $r=b$.  Then we interpolate this solution on our sub-domains to generate initial guesses for the process as we begin the method of continuation, increasing the magnitude of the angles incrementally up to $|\psi_a|,|\psi_b|$.

\begin{figure}[t]
	\centering
	\scalebox{0.35}{\includegraphics{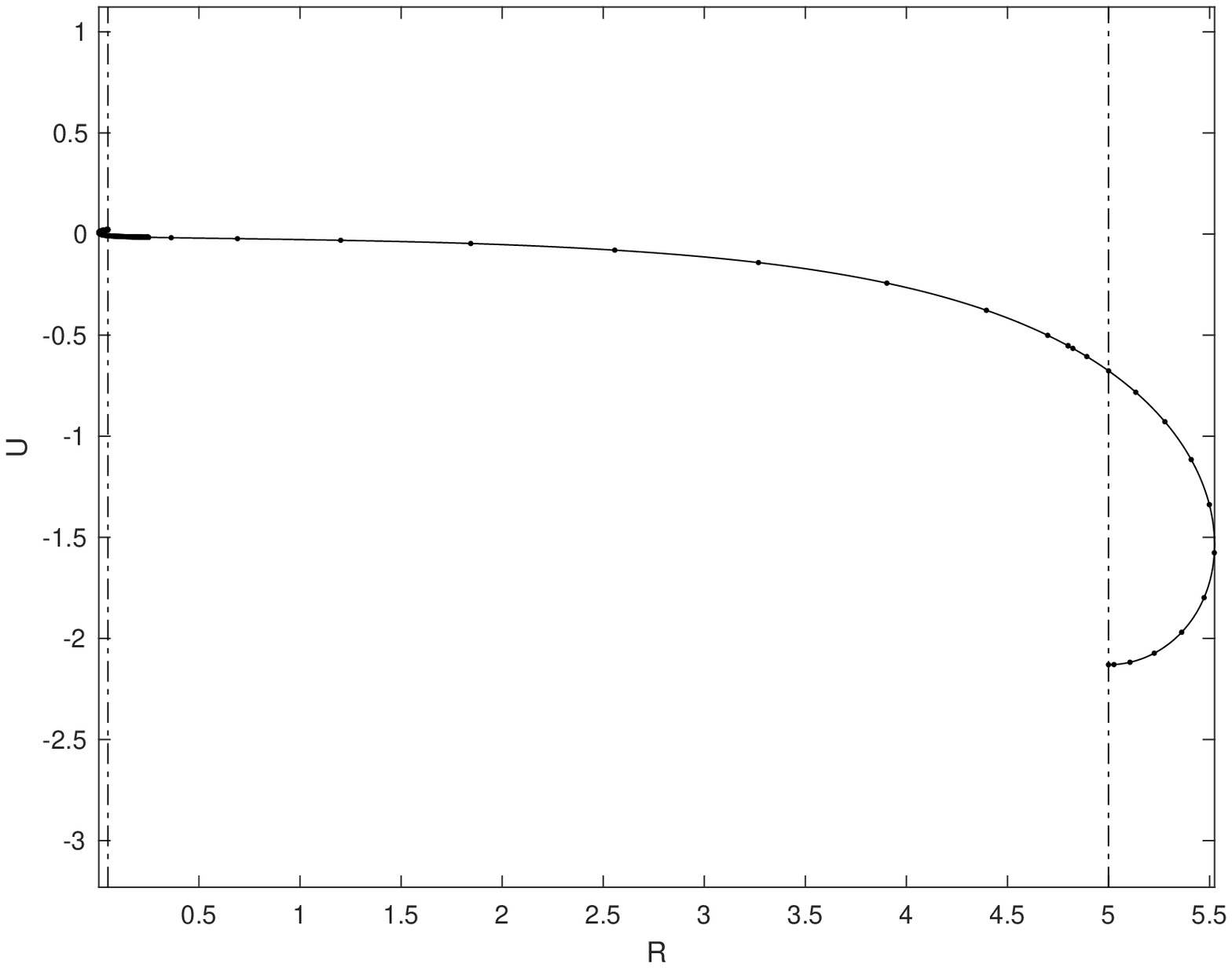}}
	\scalebox{0.35}{\includegraphics{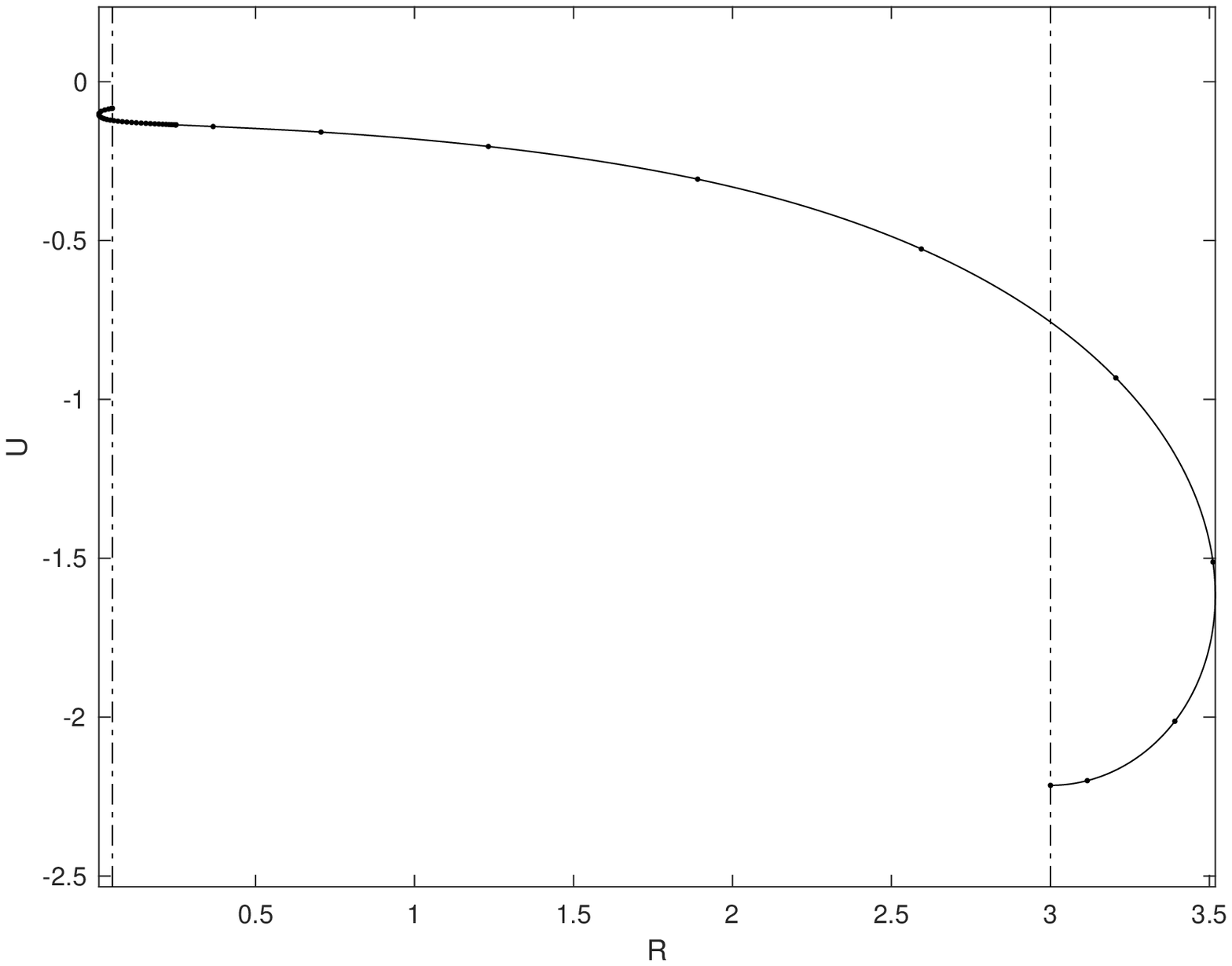}}
	\caption{A capillary surface with radii $a = 0.05$ and  $b = 5$ and inclination angles $\psi_a, \psi_b = -\pi$ using 2RZ. (L)   A capillary surface with radii $a = 0.05$ and  $b = 3$ and inclination angles $\psi_a, \psi_b = -\pi$ using A2RZ. (R)   }
	\label{fig:hook}
\end{figure}

\begin{figure}[t]
	\centering
		\scalebox{0.35}{\includegraphics{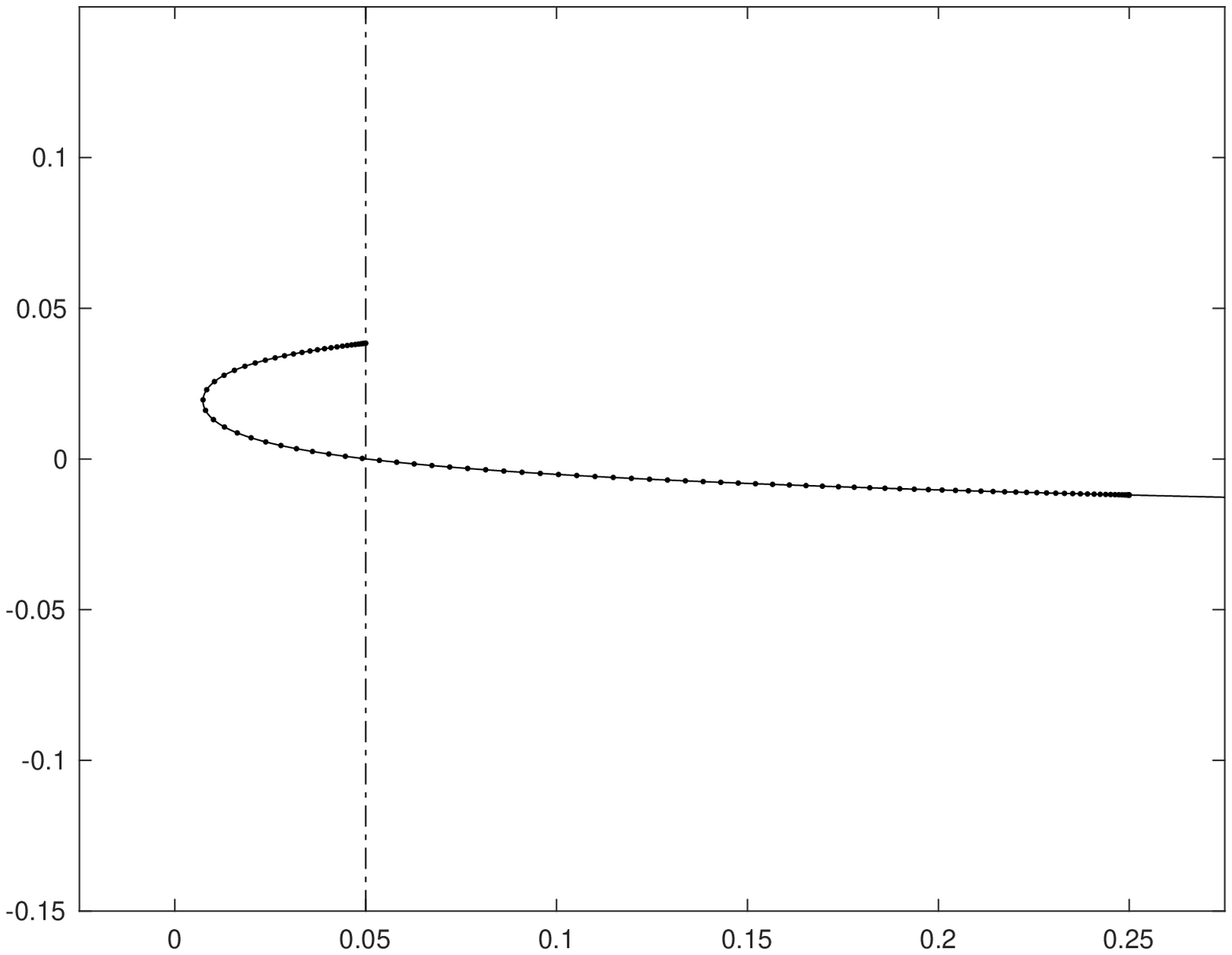}}
		\scalebox{0.35}{\includegraphics{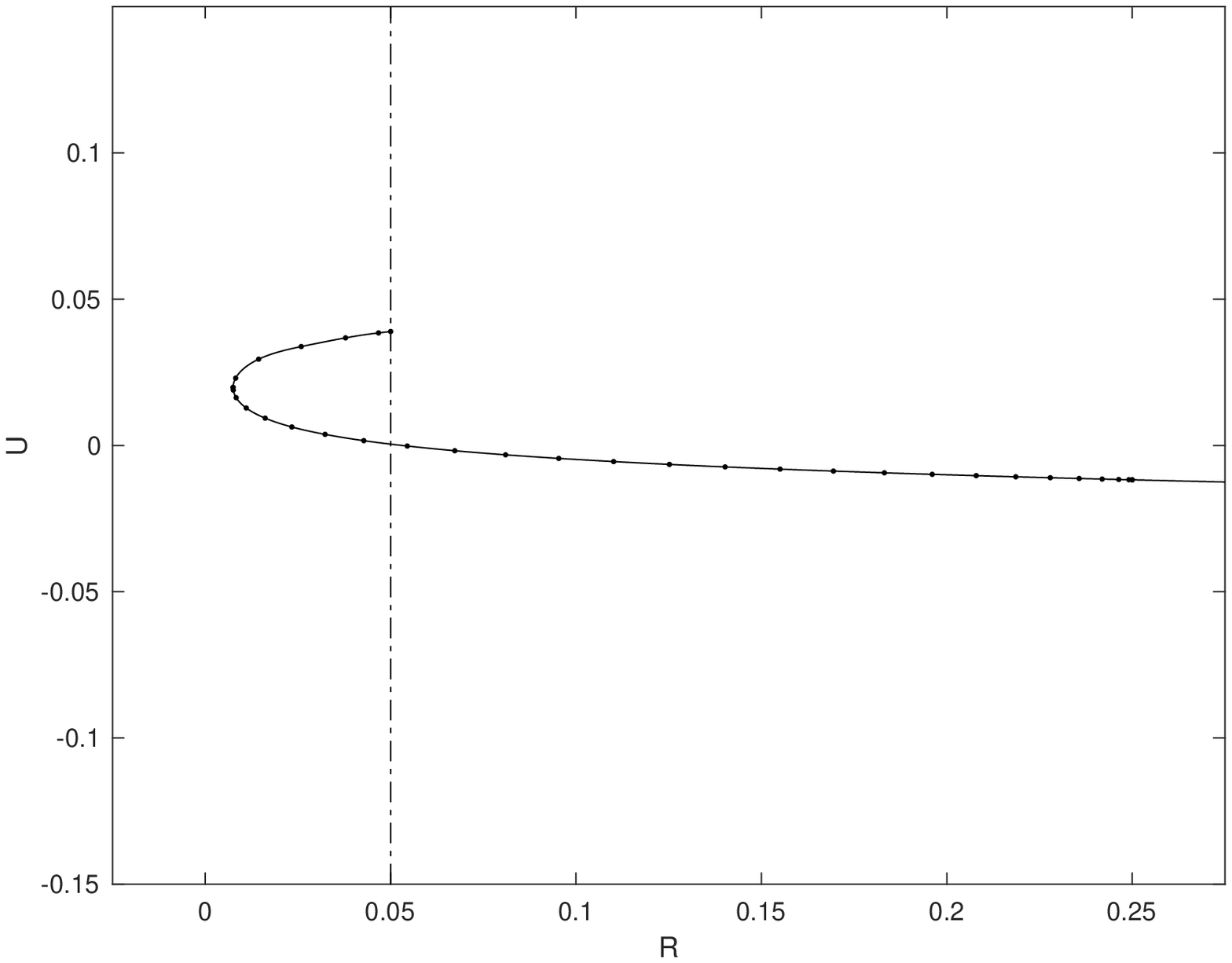}}
		\caption{Zooming in on the hook with radii $a = 0.05$ and  $b = 5$ and inclination angles $\psi_a, \psi_b = -\pi$.  The display on the left uses 2RZ and the display on the right uses A2RZ. }
		\label{fig:hookzoom}
	\end{figure}

Some examples of these experiments are collected in Tables~\ref{table2}-\ref{table3}, and can be compared with the three sub-domain results there.  Figures~\ref{fig:hook} (L) and \ref{fig:hookzoom} contain the last example from Table~\ref{table3}.  Figure~\ref{fig:hook} (R) has a  smaller value of $b$ so that the hook near $a$ is more visible.

This leads us to our final case.  Our goal here is to provide a refinement for the configuration pictured in Figure~\ref{fig:hook}.  We use the domain decomposition we just described based on two sub-domains $\Omega_1$ and $\Omega_2$ as a starting place.  Then we  mark some $\bar\psi$ and we will divide the  sub-domain $\Omega_1$ into two regions $\Omega_\alpha$ and $\Omega_\beta$, however, the starting point of these algorithms will not in general have angles that would fall in $\Omega_\alpha$.  We will proceed with the method of continuity to increase the magnitude of the prescribed angle at $r = a$ through 21 evenly spaced values starting with magnitude $\pi/2$ and ending at $|\psi_a|$,  which is presumably a value near $\pi$.  We use the two radial sub-domain formulation up to but not including the last step in this continuation.  On that last step we pick $\bar\psi = -\pi/2$ when $\psi_a < 0$, and we split $\Omega_1$ into $\Omega_\alpha$ and $\Omega_\beta$.  This gives us three sub-domains, and the boundary between $\Omega_\alpha$ and $\Omega_\beta$ is marked by $\bar\psi$, while the boundary between $\Omega_\beta$ and $\Omega_2$ is marked by $r = a + \delta$.  Then in this setting we have arc-length $\bar\ell$ for $\Omega_\alpha$ and the arc-length over $\Omega_2$ is $\ell_1 - \bar\ell$.  The matching conditions at $\bar\psi$ are given by
\begin{eqnarray}
	r_\alpha(\bar\ell) = r_\beta(\bar\ell), \\
	u_\alpha(\bar\ell) = u_\beta(\bar\ell). 
\end{eqnarray}
We then roughly proceed as before.  We first generate
\begin{equation}
	\label{NA2R}
	N(\mathbf{v}) = 
	\begin{bmatrix}
			R_\alpha^\prime(\tau_\alpha) - \ell_\alpha\cos\Psi_\alpha(\tau_\alpha)\\
			R_\beta^\prime(\tau_\beta) - \ell_\beta\cos\Psi_\beta(\tau_\beta)\\
			R_2^\prime(\tau_2) - \ell_2\cos\Psi_2(\tau_2)\\
			U_\alpha^\prime(\tau_\alpha) - \ell_\alpha\sin\Psi_\alpha(\tau_\alpha)\\
			U_\beta^\prime(\tau_\beta) - \ell_\beta\sin\Psi_\beta(\tau_\beta)\\
			U_2^\prime(\tau_2) - \ell_2\sin\Psi_2(\tau_2)\\
			\Psi_\alpha^\prime(\tau_\alpha)R_\alpha(\tau_\alpha) + \ell_\alpha\sin\Psi_\alpha(\tau_\alpha) - \kappa\ell_\alpha R_\alpha(\tau_\alpha) U_\alpha(\tau_\alpha) \\
			\Psi_\beta^\prime(\tau_\beta)R_\beta(\tau_\beta) + \ell_\beta\sin\Psi_\beta(\tau_\beta) - \kappa\ell_\beta R_\beta(\tau_\beta) U_\beta(\tau_\beta) \\
			\Psi_2^\prime(\tau_2) + \frac{\ell_2\sin\Psi_2(\tau_2)}{R_2(\tau_2)} - \kappa\ell_2 U_2(\tau_2) \\
			R_\alpha(-1) - a \\
			R_\alpha(1) - R_\beta(-1)\\
			\Psi_\alpha(-1) - \psi_a  \\
			\Psi_\alpha(1) - 	\bar\psi \\
			U_\alpha(1) - U_\beta(-1)\\
			R_\beta(1) - a - \delta\\
			\Psi_\beta(-1) - \bar\psi \\
			\Psi_\beta(1) - 	\Psi_2(-1) \\
			R_2(-1) - a - \delta\\
			U_\beta(1) - U_2(-1) \\
			R_2(1) - b \\
			\Psi_2(1) - \psi_b 
		\end{bmatrix}.
\end{equation}
Then the blocks of $L$ that need to be updated for this case are  given by
$$
M_{21}  = 
\begin{bmatrix}
	\Psi^\prime_\alpha - \kappa\ell_\alpha U_\alpha & 0& 0 & -\kappa\ell_\alpha R_\alpha & 0& 0 \\
	0 & \Psi^\prime_\beta - \kappa\ell_\beta U_\beta& 0& 0 & -\kappa\ell_\beta R_\beta & 0 \\
	0 & 0 & \frac{-\ell_2\sin\Psi_2 }{R_2^2}& 0& 0 & -\kappa\ell_2 & 0
\end{bmatrix},
$$
$$
M_{22}  = 
\begin{bmatrix}
	D  + \ell_\alpha\cos\Psi_\alpha & 0& 0&  - \kappa R_\alpha U_\alpha & 0& 0\\
	0& D  + \ell_\beta\cos\Psi_\beta& 0& 0&  - \kappa R_\beta U_\beta &  0\\
	0& 0& D  + \frac{\ell_2\cos\Psi_2 }{R_2}& 0& 0&  - \kappa U_2 
\end{bmatrix},
$$
and the blocks corresponding to boundary conditions and matching conditions are given by
$$
M_{31}  = 
\begin{bmatrix}
	D^0_{-1} & 0 & 0 & 0 & 0& 0 \\
	D^0_{1} & -D^0_{-1} & 0 & 0 & 0& 0 \\
	0 & 0 & 0 & 0 & 0& 0 \\
	0 & 0 & 0 & 0 & 0& 0 \\
	0 & 0 & 0 & D^0_1 & -D^0_{-1} & 0 \\
	0 & D^0_{1} & 0 & 0 & 0 & 0 \\
	0 & 0 & 0 & 0 & 0 & 0 \\
	0 & 0 & 0 & 0 & 0 & 0 \\
	0 & 0 & D^0_{-1} & 0 & 0 & 0  \\
	0 & 0 & 0 & 0 & D^0_1 & -D^0_{-1} \\
	0 & 0 & 0 & 0 & 0 & D^0_{1} \\
	0 & 0 & 0 & 0 & 0 & 0 
\end{bmatrix},
$$
and
$$
M_{32}  = 
\begin{bmatrix}
	0 & 0 & 0 & 0 & 0& 0 \\
	0 & 0 & 0 & 0 & 0& 0 \\
	D^0_{-1} & 0 & 0 & 0 & 0& 0 \\
	D^0_{1} & 0 & 0 & 0 & 0& 0 \\
	0 & 0 & 0 & 0 & 0& 0 \\
	0 & 0 & 0 & 0 & 0& 0 \\
	0 & D^0_{-1} & 0 & 0 & 0& 0 \\
	0 & D^0_{1} & -D^0_{-1} & 0 & 0 & 0 \\
	0 & 0 & 0 & 0 & 0& 0 \\
	0 & 0 & 0 & 0 & 0& 0 \\
	0 & 0 & 0 & 0 & 0& 0 \\
	0 & 0 & D^0_1 & 0 & 0 & 0
\end{bmatrix}.
$$
With these objects in hand, we again use the algorithm from the flowchart in Figure~\ref{flowchart}.  These results are compared to the previous approaches in Table~\ref{table3A}, and in each of the cases considered here the base code failed.

\begin{table}[t]
	\centering
	\begin{tabular}{||l|l|l|l|l|l|l|l|l|l||}
		\cline{1-10}
		$a$ & $b$ & $\psi_a$ & $\psi_b$ & base  $n_v$ & base  $n_N$ & 3RZ  $n_v$& 3RZ $n_N$& 2RZ $n_v$& 2RZ $n_N$\\
		\hline \hline
		0.05 & 1 & $-15\pi/16$ & $-15\pi/16$ & 1162 & 78   & dnc &          & 332 & 198\\
		0.05 & 1 & $-31\pi/32$ & $-\pi$             & 2314 & 82   & dnc &         & 620 & 250  \\
		0.05 & 3 & $-15\pi/16$ & $-15\pi/16$ & 2650 & 134 & 114   & 333 & 80 & 150  \\
		0.05 & 5 & $-\pi/2$        & $-\pi/2$          & 154   & 18   & 126    & 157 & 80 & 5 \\
		0.05 & 5 & $-31\pi/32$ & $-\pi$              & dnc &            & 426   & 358 & 620 & 228  \\
		\hline
	\end{tabular}
	\caption{Annular type solutions with radii $a$ and $b$, a selection of inclination angles, and the total number of points.  In this table all cases had $\delta = 0.2$.  The 3RZ cases that do not converge are due to the poor suitability of the standard initial guesses for these problems.}
	\label{table3}
\end{table}

\begin{table}[t]
	\centering
	\begin{tabular}{||l|l|l|l|l|l|l|l|l|l||}
		\cline{1-10}
		$a$ & $b$ & $\psi_a$ & $\psi_b$ &  3RZ  $n_v$& 3RZ $n_N$& 2RZ $n_v$& 2RZ $n_N$ & A2RZ  $n_v$ & A2RZ  $n_N$\\
		\hline \hline
		0.05 &	2	& $- 15\pi/16$ &	$-  15\pi/16$ 	 & 366 &	535	& 332	& 184	&	210 & 	20	\\
		0.05 & 2	& $-61\pi/64$ &  	$-\pi$ & 654 & 551	&	620	& 197	& 	354	& 25 \\		
    	0.05 & 2	& $- 31\pi/32$ & 	$-\pi$ & 654 & 577	&	620	& 227	& 	642	& 39 \\
    	\hline
		0.05 & 3	& $- 15\pi/16$ & 	$-15\pi/16$ & 114	& 359	& 80	& 150	& 	138 &	31		\\
    	0.05 & 3	& $- 61\pi/64$ & 	$-\pi$ & 366 & 404	& 	332	& 193	& 138 &	38		\\
		0.05 & 3	& $- 31\pi/32$ & 	$-\pi$ & 654 & 447	&	620	& 236	& 354	& 36	\\
		    	\hline
		0.05 & 5	& $-15\pi/16$ & 	$- 15\pi/16$ & 126	& 284	& 	80	& 145	&	138	& 29	\\
    	0.05 &	5	& $-61\pi/64$ &  	$-\pi$ & 378 &	322 & 332	& 186	&	138 & 	37		\\
		0.05 & 	5	& $-31\pi/32$ &  	$-\pi $ & 426 & 358	& 620	& 228	& 210 & 34		\\
		\hline
	\end{tabular}
	\caption{Solutions that are the image of an annulus with radii $a$ and $b$, a selection of inclination angles, and the total number of points as before.  In this table all cases had $\delta = 0.2$.}
	\label{table3A}
\end{table}

\section{Conclusions and closing remarks}
\label{conclusions}

We have shown how several domain decompositions can improve performance over the adaptive Chebyshev spectral method in \cite{Treinen2023}.  These new codes are primarily aimed at addressing problematic configurations identified in that work.  The nonlinearities in the equations lead to complicated behavior that is nontrivial to classify, and as such it is difficult to determine beforehand which approach will lead to success.  For the general scope of the annular problems, the existence theory is only partially completed, and the uniqueness of solutions is completely open \cite{Treinen2012}.   With this perspective, it may take some experimenting with the number of sub-domains and the location of the corresponding boundary points to achieve success.

\begin{figure}[t]
	\centering
	\scalebox{0.35}{\includegraphics{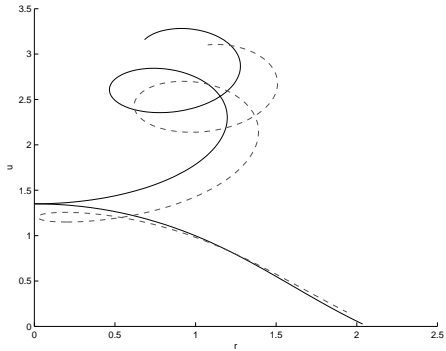}}
		\caption{Two curves from the families of global solutions of \eqref{drds}-\eqref{dpsids}.  One curve has a value of $r = 0$ with two curves horizontally tangent there, and the other has a minimum value of $r$ that corresponds to a loop containing a vertical point. }
		\label{fig:axis_curves}
	\end{figure}

For example, in the configurations where $\psi_a, \psi_b =-\pi$ when  $0<a\ll 1$ and $b$ is moderately sized the performance of the base code \cite{Treinen2023}   is poor or fails completely, depending on how small $a$ is.  These configurations correspond to a singular limiting process that was explored by Bagley and Treinen \cite{BagleyTreinen2018}.  In that work the global solutions of \eqref{drds}-\eqref{dpsids} were considered for all arc-lengths.  These global solutions were classified into families of solutions, and the leftmost vertical point of the interface was a parameter used in this classification.  When the radius of that leftmost vertical point is positive the curve forms a loop there.  As the radius of that vertical point goes to zero the loop collapses into two curves that are tangent at $r=0$, as indicated in Figure~\ref{fig:axis_curves}.   (A rigorous study of this phenomenon has not yet been undertaken.)   This behavior highlights the difficulty of the problem.  While our attempts to use this multi-scale approach on this problem has improved the performance, we are still unable to find solutions up to $\psi_a = -\pi$ in some cases.

Our 3RZ code is designed to address $a\ll \infty$ and $b$ somewhat large.  Table~\ref{table2} shows that this is working well.  It does not function particularly well when $b$ is closer to $a$, as seen in Table~\ref{table3}.  This could be improved by using the initial guess scheme introduced for the 2RZ case.  

Our 2RZ code is designed to address $0<a\ll 1$ with moderate values of $b$, and while it does well with this difficult problem, getting to $\psi_a = -\pi$ is quite difficult.  The A2RZ case is an attempt to resolve this, and while there is further improvement, getting all the way to $\psi_a = -\pi$ in all cases is still illusive.

This paper has focused on certain choices of inclination angles at particular radii.  The opposite choices for those inclination angles lead to reflections of the surface over the radial axis, and work equally well.  In general there are four parameters to explore, and it is possible that there are other difficult cases that we have not yet encountered.
The programs we have written are available at the GitHub site referenced in the introduction and could be modified to treat these new cases if they are discovered.


\end{document}